# Model selection and sensitivity analysis for sequence pattern models[*]

## Mayetri Gupta[1]

*University of North Carolina at Chapel Hill*

**Abstract:** In this article we propose a maximal a posteriori (MAP) criterion for model selection in the motif discovery problem and investigate conditions under which the MAP asymptotically gives a correct prediction of model size. We also investigate robustness of the MAP to prior specification and provide guidelines for choosing prior hyper-parameters for motif models based on sensitivity considerations.

## 1. Introduction: statistical challenges in motif discovery

Genome sequencing projects have led to a rapid growth of databases of genome sequences for DNA, RNA and proteins. The task of extracting insight into gene regulatory networks from these massive amounts of data represents a major scientific challenge. A first step towards understanding the process of gene regulation is the identification of short recurring patterns (about 8-20 nucleotides long), called motifs, in a set of bio-polymer sequences. Motifs correspond to functionally important parts of molecules, such as the sites where certain proteins, called transcription factors (TFs) bind, to control gene expression. In spite of a plethora of computational methods developed to address the problem of motif discovery (Sandve and Drablos [15] and Gupta and Liu [7]), most approaches still suffer from a lack of specificity in motif predictions- yielding a high number of false positives that appear to have no known biological significance. Thus, one of the fundamental questions that arise is whether the patterns discovered from the sequence data by an algorithm are real. Although confirming the biological relevance of these findings often requires further biological experimentation, it is at least important to assess their statistical significance. We approach this question of model selection from a Bayesian viewpoint, using an analytical approximation to the Bayes factor– the maximal a posteriori (MAP) scoring criterion. As an evaluation of its performance, it is shown that the MAP asymptotically attains several desirable properties. Since the Bayes factor or any such Bayesian model selection criterion necessarily involves parameters of the prior distribution, we also conduct sensitivity analyses to judge the effect of prior parameter specifications on posterior inference and prescribe robust hyperparameter choices for practitioners.

[*]Supported by an IBM junior faculty award from the University of North Carolina at Chapel Hill.
[1]Department of Biostatistics, 3107C McGavran Greenberg CB# 7420, Chapel Hill, NC 27599-7420, USA, e-mail: gupta@bios.unc.edu

### 1.1. Bayesian stochastic dictionary model for sequence patterns

For convenience, we view the sequence data as a *single* sequence $\mathcal{S} = \{x_1 x_2 \ldots x_n\}$ of length $n$. $\mathcal{S}$ is assumed to be generated by the concatenation of words from a *dictionary* $\mathcal{D}$ of size D, where $\mathcal{D} = \{M_1, M_2, \ldots, M_D\}$, sampled randomly according to a probability vector $\boldsymbol{\rho} = (\rho(M_1), \ldots, \rho(M_D))$. The likelihood of $\mathcal{S}$ is

$$(1.1) \qquad P(\mathcal{S} \mid \boldsymbol{\rho}) = \sum_{\Pi} \prod_{i=1}^{N(\Pi)} \rho(\mathcal{S}[P_i]) = \sum_{\Pi} \prod_{j=1}^{D} \rho(M_j)^{N_{M_j}(\Pi)},$$

where $\Pi = (P_1, \ldots, P_k)$ is a partition of the sequence so that each part $P_i$ corresponds to a word in the dictionary, $N(\Pi)$ is the total number of partitions in $\Pi$, and $N_{M_j}(\Pi)$ is the number of occurrences of word $M_j$ in the partition.

Assume that the first $b$ ($b < D$) words in the dictionary are the single letters (for DNA, $b = 4$ letters, $\{A, C, G, T\}$). Let $\boldsymbol{\rho} = (\rho_1, \ldots, \rho_D)$ denote the word usage probabilities for the set of all words in the dictionary. If the partition $\Pi = (P_1, \ldots, P_k)$ of the sequence into words were known, the resulting distribution of counts of words $\boldsymbol{N} = (N_1, \ldots N_D)^T$ would be multinomial characterized by the probability vector $\boldsymbol{\rho}$. Suppose we have $D - b$ motifs (words other than single letters) of widths $w_{b+1}, \ldots, w_D$. In the *stochastic dictionary* framework, "words" $M_j$, $(j = b+1, \ldots, D)$ are stochastic matrices denoted by $\{\boldsymbol{\Theta}_{b+1}, \ldots, \boldsymbol{\Theta}_D\} = \boldsymbol{\Theta}^{(D)}$. For the $k^{\text{th}}$ word of width $w_k$, its probability matrix is represented as $\boldsymbol{\Theta}_k = (\boldsymbol{\theta}_{1k}, \ldots, \boldsymbol{\theta}_{w_k k})$, each of the $w_k$ columns giving the probabilities of occurrence of the four letters at that position of the word. When multiple occurrences of word $k$ are aligned, the letter counts in the $j^{\text{th}}$ aligned column, $\mathbf{c}_{jk} = (c_{1jk}, \ldots c_{bjk})^T$, $(j = 1, \ldots w_k)$, are characterized by the probability vectors $\boldsymbol{\theta}_{jk} = (\theta_{1jk}, \ldots, \theta_{bjk})^T$ $(j = 1, \ldots w_k)$, of a *product* multinomial model. The count matrices corresponding to the motifs are denoted as $\{\boldsymbol{C}_{b+1}, \ldots, \boldsymbol{C}_D\} = \mathcal{C}$.

The partition variable $\Pi$ can be equivalently expressed by the motif site indicators, denoted as $\boldsymbol{A} = \{A_{ik}; i = 1, \ldots, n, k = b+1, \ldots D\}$, where $A_{ik} = 1$ (0) if $i$ is the start of a site corresponding to motif type $k$ (otherwise). The complete data likelihood then is: $L(\boldsymbol{N}, \mathcal{C}, \boldsymbol{A} \mid \boldsymbol{\Theta}^{(D)}, \boldsymbol{\rho}) \propto \prod_{l=1}^{D} \rho_l^{N_l} \prod_{k=b+1}^{D} \prod_{j=1}^{w_k} \prod_{i=1}^{b} \theta_{ijk}^{c_{ijk}}$. We assume a Dirichlet prior distribution for $\boldsymbol{\rho}$, $\boldsymbol{\rho} \sim \text{Dirichlet}(\boldsymbol{\beta}_0)$, $\boldsymbol{\beta}_0 = (\beta_{01}, \ldots \beta_{0D})$, and a corresponding product Dirichlet prior (i.e., independent priors over the columns) PD($\boldsymbol{B}$) for $\boldsymbol{\Theta}_k = (\boldsymbol{\theta}_{1k}, \ldots, \boldsymbol{\theta}_{w_k k})$, $(k = b+1, \ldots D)$. A Dirichlet prior is a natural choice for this application, modeling the joint prior densities of proportions, and is conjugate with the likelihood, leading to easy computation of the posterior densities. Let $\boldsymbol{B} = (\boldsymbol{\beta}_1, \boldsymbol{\beta}_2, \ldots \boldsymbol{\beta}_{w_k})$ be a $b \times w_k$ matrix with $\boldsymbol{\beta}_j = (\beta_{1j}, \ldots \beta_{bj})^T$. The conditional posterior distribution of $\boldsymbol{\Theta}_k \mid \boldsymbol{N}, \mathcal{C} \propto \prod_{j=1}^{w_k} \prod_{i=1}^{b} \theta_{ijk}^{c_{ijk} + \beta_{ij}}$, which is product Dirichlet PD($\boldsymbol{B} + \boldsymbol{C}_k$), the pseudo-count parameters $\boldsymbol{B}$ being updated by the column counts of the $k^{\text{th}}$ word, $\boldsymbol{C}_k = (\mathbf{c}_{1k}, \ldots \mathbf{c}_{w_k k})$. The conditional posterior of $\boldsymbol{\rho} \mid \boldsymbol{N}, \mathcal{C}$ is Dirichlet($\boldsymbol{N} + \boldsymbol{\beta}_0$) $\propto \prod_{l=1}^{D} \rho_l^{N_l + \beta_{0l}}$. Bayes estimates of $\mathcal{O} = (\boldsymbol{\Theta}, \boldsymbol{\rho})$ may be obtained through a DA procedure using the conditional distributions $P(\boldsymbol{A} \mid \mathcal{O}, \mathcal{S})$ and $P(\mathcal{O} \mid \mathcal{S}, \boldsymbol{A})$, using techniques of dynamic programming (DP) (Gupta and Liu [6]).

## 2. Model selection in the motif discovery problem

A traditional model selection approach is to use likelihood-based criteria, for instance, penalized likelihoods. Two widely used criteria are the (i) AIC (Akaike's

Information Criterion) and (ii) BIC (Bayesian Information or Schwarz Criterion). Leroux [12] proved that, under mild conditions, the estimator obtained with the number of components selected using AIC or BIC in mixture models is consistent for the true mixing distribution, and has in the limit at least as many components. In particular examples, however, it has been observed that the AIC and BIC do not give identical results: in our case, the BIC may heavily penalize longer motifs, when the data set size is large. An alternative for judging the degree of conservation of a motif is the Kullback-Leibler information criterion (KLI), that measures the "entropy" distance of the motif from the background: $KLI = \sum_{i=1}^{w} \sum_{j=1}^{q} \hat{\theta}_{ij} \log \frac{\hat{\theta}_{ij}}{\theta_{0j}}$. As the exact likelihood can be calculated computationally, using the procedure mentioned in Gupta and Liu [6], it is possible to calculate any of the above criteria. However, the dependence structure in the sequence models, when the position of motif sites is unknown, makes it difficult to accurately determine the sample size $n$ for BIC. Also, the KLI is sensitive to the width of the motif- making it difficult to compare the values across different motif widths. The relative performance of several of these model selection criteria on experimentally verified models for real data sets are given in Section 2.4.

## 2.1. The Bayesian approach

We can alternatively formulate the question as a Bayesian model selection problem. In the simplest scenario, it is of interest to assess whether the sequence data should be explained by model $\mathcal{M}_1$, which assumes the existence of a single nontrivial motif, or $\mathcal{M}_0$, which says that the sequences are generated entirely from a background model (e.g., an i.i.d. or Markov model). The Bayes factor, which is the ratio of the marginal likelihoods under the two models, can be computed as

$$
\begin{aligned}
\frac{p(\boldsymbol{S} \mid \mathcal{M}_1)}{p(\boldsymbol{S} \mid \mathcal{M}_0)} &= \frac{\sum_{\boldsymbol{A}} \int_{\mathcal{O}} p(\boldsymbol{A}, \boldsymbol{S}, \mathcal{O} \mid \mathcal{M}_1) d\mathcal{O}}{\int_{\mathcal{O}} p(\boldsymbol{S}, \mathcal{O} \mid \mathcal{M}_0) d\mathcal{O}} \\
&= \frac{\sum_{\boldsymbol{A}} p(\boldsymbol{A}, \boldsymbol{S} \mid \mathcal{M}_1)}{p(\boldsymbol{S} \mid \mathcal{M}_0)} \left( \equiv \frac{c_1}{c_0} \right).
\end{aligned}
$$
(2.1)

The individual additive terms in the numerator, after integrating out $\mathcal{O}$, consist of ratios of products of gamma functions. To evaluate this sum exhaustively over all partitions involves prohibitive amounts of computation. We thus need to resort to either Monte Carlo or some approximations.

It can be observed that $p(\boldsymbol{S} \mid \mathcal{M}_1)$ is the normalizing constant of $p(\boldsymbol{A} \mid \boldsymbol{S}, \mathcal{M}_1) = \frac{p(\boldsymbol{A}, \boldsymbol{S} \mid \mathcal{M}_1)}{p(\boldsymbol{S} \mid \mathcal{M}_1)} = c_1 q(\boldsymbol{A} \mid \boldsymbol{S}, \mathcal{M}_1)$, (say) where $p(\boldsymbol{A}, \boldsymbol{S} \mid \mathcal{M}_1)$ is of known form. Computational approximations to estimate normalizing constants is a subject of much recent research e.g. Meng and Wong [13], Chen and Shao [3, 4], Chib and Jeliazkov [5]. We tested importance and bridge sampling (Meng and Wong [13]) methods for estimating the ratio, which are possible to implement since we can obtain MCMC draws from $p_1(\boldsymbol{A} \mid \boldsymbol{S})$. However, neither method performs well in most of the applications, probably since the MCMC chain is extremely sticky and hence the draws do not well represent the true distribution. The high concentration of the density of $\boldsymbol{A}$ around the mode motivates us to seek a measure of significance that uses the modal information, and such a method is next discussed.

## 2.2. The MAP approximation to the Bayes factor

An obvious lower bound for (2.2) is $p(\boldsymbol{A}^*, \boldsymbol{S} \mid \mathcal{M}_1)/p(\boldsymbol{S} \mid \mathcal{M}_0)$, where $\boldsymbol{A}^*$ is the maximizer of the ratio. We demonstrate how this lower bound, called the *maximal a posteriori* score (denoted by MAP($\boldsymbol{A}^*$)), can be used as a model selection criterion.

**Definition 1.** The Maximal a Posteriori (MAP) score under the stochastic dictionary model $\mathcal{M}_1$ with a single motif of width $w$ and $d$ letters at the "best" alignment $A^*$, compared to the background model $\mathcal{M}_0$ (with alphabet size $d$) is given by

$$(2.2) \qquad MAP(A^*) = \frac{P(\boldsymbol{S}, A^* \mid \mathcal{M}_1)}{P(\boldsymbol{S} \mid \mathcal{M}_0)} = \frac{p_1(\boldsymbol{S}, A^*)}{p_0(\boldsymbol{S})}.$$

Let $\mathbf{c} = \sum_{j=1}^{w} \mathbf{c}_j$ denote the word matrix column counts, $\boldsymbol{N}_{[1:d]} = \{N_1, \ldots, N_d\}$ the letter frequencies and $\boldsymbol{\beta}_{0[1:d]} = \{\beta_{01} \ldots, \beta_{0d}\}$ the prior pseudo-counts for the $d$ letters, $\boldsymbol{\beta}_0$ denoting the pseudo-counts for the entire set of $D$ words ($D = d+1$ for 1 long word). Then
(2.3)
$$MAP(A^*) = \frac{\int_{\Theta, \boldsymbol{\rho}_1} \prod_{i=1}^{d+1} \rho_{1i}^{N_i(A^*, \mathcal{M}_1)} \prod_{j=1}^{w} \prod_{k=1}^{d} \Theta_{jk}^{C_{jk}(A^*, \mathcal{M}_1)} p(\boldsymbol{\rho}_1, \boldsymbol{\Theta}) d\boldsymbol{\rho}_1 d\Theta}{\int_{\boldsymbol{\rho}_0} \prod_{i=1}^{d} \rho_{0i}^{N_i(\mathcal{M}_0)} p(\boldsymbol{\rho}_0) d\boldsymbol{\rho}_0},$$

where $N(\cdot)$ and $C(\cdot)$ denote the word counts and motif column counts respectively, as defined in Section 1.1. For any indicator vector $\boldsymbol{A}$ (representing the sampled motif locations), we can integrate out $\mathcal{O}$ and obtain

$$\begin{aligned} \log \frac{p(\boldsymbol{A}, \boldsymbol{S} \mid \mathcal{M}_1)}{p(\boldsymbol{S} \mid \mathcal{M}_0)} &= \log \left\{ \frac{\Gamma(\boldsymbol{N} + \boldsymbol{\beta}_0)}{\Gamma(|\boldsymbol{N} + \boldsymbol{\beta}_0|)} \frac{\Gamma(|\boldsymbol{\beta}_0|)}{\Gamma(\boldsymbol{\beta}_0)} \right\} \\ &\quad - \log \left\{ \frac{\Gamma(\boldsymbol{N}_{[1:d]} + \mathbf{c} + \boldsymbol{\beta}_{0[1:d]})}{\Gamma(|\boldsymbol{N}_{[1:d]} + \mathbf{c} + \boldsymbol{\beta}_{0[1:d]}|)} \frac{\Gamma(|\boldsymbol{\beta}_{0[1:d]}|)}{\Gamma(\boldsymbol{\beta}_{0[1:d]})} \right\} \\ &\quad + \sum_{j=1}^{w} \log \left\{ \frac{\Gamma(\mathbf{c}_j + \boldsymbol{\gamma})}{\Gamma(|\mathbf{c}_j + \boldsymbol{\gamma}|)} \frac{\Gamma(|\boldsymbol{\gamma}|)}{\Gamma(\boldsymbol{\gamma})} \right\}. \end{aligned}$$

When there are $D - d > 1$ motifs (nontrivial stochastic words) in the dictionary, the logarithm of the MAP score (denoted logMAP) can be computed as:

$$\begin{aligned} \mathrm{logMAP}(\boldsymbol{A}) &= \log \left\{ \frac{\Gamma(\boldsymbol{N} + \boldsymbol{\beta}_0)}{\Gamma(|\boldsymbol{N} + \boldsymbol{\beta}_0|)} \frac{\Gamma(|\boldsymbol{\beta}_0|)}{\Gamma(\boldsymbol{\beta}_0)} \right\} \\ &\quad - \log \left\{ \frac{\Gamma(\boldsymbol{N}_{[1:d]} + \mathbf{c} + \boldsymbol{\beta}_{0[1:d]})}{\Gamma(|\boldsymbol{N}_{[1:d]} + \mathbf{c} + \boldsymbol{\beta}_{0[1:d]}|)} \frac{\Gamma(|\boldsymbol{\beta}_{0[1:d]}|)}{\Gamma(\boldsymbol{\beta}_{0[1:d]})} \right\} \\ &\quad + \sum_{k=b+1}^{D} \sum_{j=1}^{w_k} \log \left\{ \frac{\Gamma(\mathbf{c}_{jk} + \boldsymbol{\gamma})}{\Gamma(|\mathbf{c}_{jk} + \boldsymbol{\gamma}|)} \right\} + \sum_{k=b+1}^{D} w_k \log \left\{ \frac{\Gamma(|\boldsymbol{\gamma}|)}{\Gamma(\boldsymbol{\gamma})} \right\}. \end{aligned}$$

where now $\mathbf{c} = \sum_{k=d+1}^{D} \sum_{j=1}^{w} \mathbf{c}_{jk}$, and $\boldsymbol{N}_{[1:d]}$ and $\beta_{0[1:d]}$ are the same as above.

In most simulations and actual data sets, we observed that (i) the MAP score dominates the other components of the Bayes factor and is often more reliable to differentiate the two models than the approximated Bayes factor via importance sampling and (ii) the observed logMAP score in i.i.d. data sets (containing no "true" motifs) for any selected alignment was always lower than that of the null alignment.

## 2.3. Asymptotic results for MAP scoring

We now derive asymptotic optimality criteria for judging the performance of the MAP scoring statistic as a model selection criterion. In this section, we assume that $\boldsymbol{\rho}_0 \mid \mathcal{M}_0 \sim \text{Dir}(\boldsymbol{\alpha}^{(d)})$ and $\boldsymbol{\rho}_1 \mid \mathcal{M}_1 \sim \text{Dir}(\boldsymbol{\beta}^{(d+1)})$, where $\boldsymbol{\alpha} = (\alpha_1, \ldots \alpha_d)^T$ and $\boldsymbol{\beta} = (\beta_1, \ldots \beta_{d+1})^T$. Additionally, assume that the sequence is generated by *random* draws from an alphabet of size $d$ according to some arbitrary probabilities. Let us denote the score for the "true" alignment as $MAP(A^0)$ (so that $MAP(A^*) \geq MAP(A^0)$). The key result for establishing the asymptotic optimality of the MAP criterion for model selection is as follows.

**Theorem 2.1.** *For $m$ exact motif instances, each of length $w$, in a dataset of size $N$, let $\frac{m}{N} \xrightarrow{P} c$, $(c < \frac{1}{w})$, $\frac{N_{0i} - N_{1i}}{m} \xrightarrow{P} k_i w$, and $\frac{N_{0i}}{N} \xrightarrow{P} \theta_{0i}$, (where $N_{0i}$ and $N_{1i}$ are the frequency of letter $i$ in the background and motif respectively, $\sum_{i=1}^{d} k_i = 1$, and $\sum_{i=1}^{d} \theta_{0i} = 1$.) In other words, the proportion of motif sequence to total sequence data tends to a constant; the proportion of letter $i$ in a motif tends to a constant $k_i$ and the sample proportions of letters approaches the population proportion. Also, let*

$$
\begin{aligned}
r &= c \log c + \sum_{i=1}^{d} (\theta_{0i} - k_i c w) \log(\theta_{0i} - k_i c w) \\
&\quad - \sum_{i=1}^{d} \theta_{0i} \log \theta_{0i} - [1 - c(w-1)] \log[1 - c(w-1)].
\end{aligned}
\tag{2.4}
$$

*Then, $r > 0$ is a necessary and sufficient condition for $MAP(A^0) \xrightarrow{P} \infty$ (as $N \to \infty$, $m \to \infty$) at a rate exponential in $rN$.*

If the conditions of Theorem 2.1 are satisfied, $MAP(A^*)$, and thus also the Bayes factor, approaches infinity at an exponential rate– indicating that asymptotically the MAP gives the correct model judgment. Following from this result, we denote the expression $r$ in Theorem 2.1 as the MAP divergence factor, or $MAP_{DF}$.

**Theorem 2.2.** *The maximum value attained by $r$ in expression (2.4) of Theorem 2.1 for fixed $c$, $d$ and $w$ is*

$$c \log c + (1 - cw) \log(1 - cw) + cw \log d - [1 - c(w-1)] \log[1 - c(w-1)],$$

*attained for $\theta_{0i} = k_i = \frac{1}{d}$ for $i = 1, \ldots, d$.*

It is interesting to note some special cases where Theorem 2.2 holds.

**Case 1: Symmetric base composition, symmetric motif contribution**
With equal background proportions of letters, i.e. $\theta_{0i} = \frac{1}{d}$ for $i = 1, \ldots, d$, and an equal contribution of each base to motif $k_i = \frac{1}{d}$ for $i = 1, \ldots, d$, the conditions of Theorem 2.2 are satisfied. The positivity condition of expression (2.4) reduces to

$$
\begin{aligned}
&c \log c + (1 - cw) \log(1 - cw) \\
&\quad + cw \log d - [1 - c(w-1)] \log[1 - c(w-1)] > 0.
\end{aligned}
\tag{2.5}
$$

From Theorem 2.2 it is clear that the MAP score for a motif with other than symmetric composition will diverge more slowly. Another extreme case is a "repeat-based" motif, composed of a single letter repeated throughout the word.

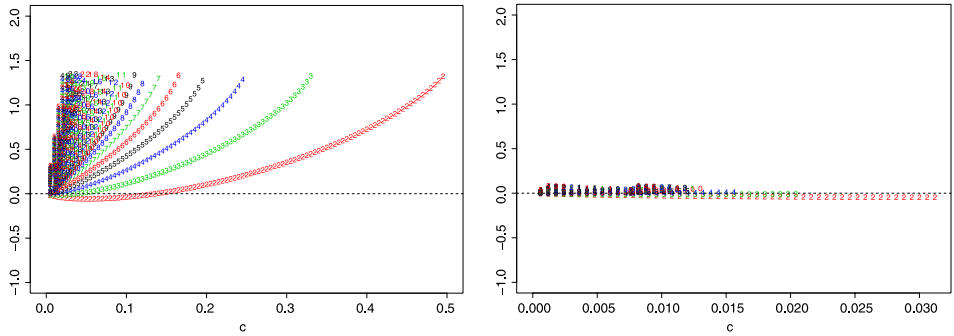

FIG 1. *(a) Values of MAP divergence factor ($MAP_{DF}$) increase with motif widths* 2 *to* 50 *(numbers on the figure) and values of c for "equal contribution" motif, i.e. an* exact *motif in which each of the d letters occurs the same number of times. (b) Values of $MAP_{DF}$ for motif widths* 2 *to* 50 *and possible values of c for a "repeat-based" motif. Repeat motifs score lower than motifs with equal base composition.*

**Case 2: Symmetric base composition, "repeat" motif**
With equal background proportions of letters, i.e. $\theta_{0i} = \frac{1}{d}$ for $i = 1, \ldots, d$, and a *repeat* motif: e.g., $k_1 = d$; $k_i = 0$ for $i > 1$, the condition in (2.4) reduces to

$$c \log c + \frac{1}{d} \log d + \left(\frac{1}{d} - dcw\right) \log \left(\frac{1}{d} - dcw\right)$$
$$-[1 - c(w-1)] \log[1 - c(w-1)] > 0.$$

(In this case, $cw$ is necessarily $< \frac{1}{d^2}$.) It is interesting to observe that a repeat motif will score *lower* than a motif pattern with a more varied representation of letters (Panel (b) of Figure 1) when all other parameters remain the same.

Next, we check how these conditions perform in practice for given $c$ and $w$. More precisely, we check the positiveness of the MAP divergence factor ($MAP_{DF}$). Except for very low widths ($w < 5$), the $MAP_{DF}$ is almost always positive (indicated in Figure 1), suggesting that the performance of the MAP score as a model selection criterion should improve with increasing motif width, which is also observed in empirical studies. Figure 1 shows a comparatively lower value of $MAP_{DF}$ for each $w$, and a slower increase with $w$ for a repeat-based motif. In both cases $MAP_{DF}$ increases with increase in motif proportion $c$, though the increase is far more marked for the motif with equal base composition.

*2.3.1. Multiple motif MAP monotonic property*

In the progressive updating method that is used to include new motifs into the dictionary, a new motif is judged to be "significant" based on the increase in the MAP score after including the new motif. In this section we show that asymptotically, the MAP score monotonically increases with the size of the true dictionary, hence the MAP *difference* criterion theoretically tends to include "true" motifs.

Let us define $MAP_q^k$ as the MAP score corresponding to the "true" alignment for $q$ word types (including $d$ single letters and $q - d$ motifs), where the $q^{\text{th}}$ word has $k$ instances in the data set. If $\beta_i > 1 \ \forall i = 1, \ldots, d$, it is trivial to show that

$$MAP_q^0 < MAP_{q-1}^l \quad \forall l > 0,$$

which means that if no motifs of the new type $q$ are present, the MAP score should theoretically show a decrease. Next, we derive conditions very similar to (2.1) which show that if the word type $q$ is truly present, the MAP score should increase– specifically, we show that $\log MAP_q^l - \log MAP_{q-1}^l \xrightarrow{P} \infty$.

**Theorem 2.3.** *Let* $\frac{m}{N} \xrightarrow{P} \rho_{q+1}$, $\frac{N_{1i}}{N} \xrightarrow{P} \rho_i$ *(for* $i = d+1, \ldots, q$*) where* $N = \sum_{i=1}^{d} N_{1i}$. *Then* $MAP_q^m / MAP_{q-1}^l \xrightarrow{P} \infty$ *(as* $m, N \to \infty$*) at an exponential rate in* $rN$ *iff* $r > 0$, *where* $r$ *is given by:*

$$\sum_{i=1}^{d} \{(\rho_i - \kappa_i w \rho_{q+1}) \log(\rho_i - \kappa_i w \rho_{q+1}) - \rho_i \log \rho_i\}$$
$$+ \log \rho_{q+1}^{\rho_{q+1}} + \log[1 - (w-1)\rho_{q+1}]^{[1-(w-1)\rho_{q+1}]}.$$

Proofs of Theorems 2.1 and 2.3 are given in the Appendix.

### 2.4. Performance of model selection criteria

For empirical comparisons, we used four data sets: CRP (S. cerevisiae or yeast), GAL4 (yeast), $\sigma_A$ (*B. subtilis*), and skeletal muscle TF MEF2 (human) (Wasserman et al. [17]). The value of the MAP divergence factor ($MAP_{DF}$) in Theorem 2.1 is shown in Figure 2, for motifs in the 4 data sets. The $MAP_{DF}$ is calculated for motif composition frequency $\mathbf{k} = (k_1, \ldots, k_b)$ matching that of the consensus motif (Columns 5 to 8 in Table 1). For all motif occurrence proportions $c$, it appears that the $MAP_{DF}$ exceeds zero for the 3 data sets with longer motifs, viz. GAL4, CRP and MEF2– which gives an indication that the consensus motif is likely to be evaluated significant (under the conditions of Theorem 2.1). For $\sigma_A$ (with lowest motif width of 6), the value of $MAP_{DF}$ is slightly below zero for very small $c$, indicating that finding the true motif may be more difficult.

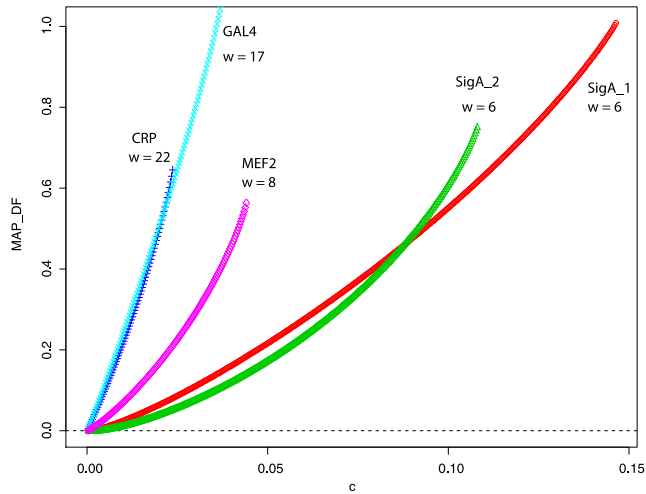

FIG 2. *Values of $MAP_{DF}$ against motif proportion $c$ for 4 data sets: (1) CRP (2) GAL4 (3) 2 motifs corresponding to $\sigma_A$ binding sites in* B. subtilis *and (4) MEF2 in human sequences. w denotes the likely width of the motif based on experimental data. ($\sigma_A^{(1)}$ denotes the motif "TATAAT" and $\sigma_A^{(2)}$ denotes "TTGACA".)*

The MAP was also compared to several likelihood-based criteria in terms of model selection performance by using the knowledge of experimentally detected motifs for these data sets. If we set our stopping decision rule to be selection of motifs that keep the logMAP score positive, this contains the set of all true motifs, for every data set, with an increased power of discrimination with increasing motif length, which agrees with our theoretically obtained results. The performance of the BIC drastically worsens when the motif length increases- perhaps because of the higher penalty for the increased number of parameters in the model.

An objection that may be raised on using the MAP for model selection is that it involves the evaluation of a "modal" value instead of a more typical Bayesian quantity involving a sample from the population. A sample-based statistic is often desirable in order to incorporate the posterior variability of the distribution into the measure used. There are two main reasons why we prefer the use of the MAP score. First, the analytical intractability of the Bayes factor forces us to choose an alternative approach – either computational or analytical. Second, it is difficult to get a "good" representative posterior sample from the distribution of interest i.e. $p(\boldsymbol{A} \mid \boldsymbol{S}, \mathcal{O}_j, \mathcal{M}_1)$, hence the lack of a "good" mean estimate. The contribution of the MAP to the BF increases as the likelihood surface becomes more spiky (hence resulting in the MAP score providing a better approximation) and simultaneously, the computational estimates for the BF are likely to perform worse, as it is more difficult to obtain a true representative sample from the distribution.

## 3. Sensitivity analysis for prior influence on the MAP criterion

In a Bayesian analysis, it is of concern to check the robustness of posterior inferences to the prior specification. Even asymptotically, Bayes factors remain sensitive to the choice of prior, even though posterior means may not (Kass [8]). Conjugate priors are not automatically robust as the tails are typically of the same form as the likelihood function, and the prior remains influential when the likelihood function is concentrated in the prior tail (Berger [2]). With a large dimensional parameter space, we need to rely on the conjugate form of the prior for analytical and computational tractability. However, it is still of interest to see the influence of the choice of prior parameters within this class on the posterior inferences, and develop an idea of situations leading to highest posterior robustness. For investigating posterior robustness over a class of priors, it is attractive to consider the $\epsilon$-contamination class (Berger [2]), defined by, $\Gamma = \{\pi : \pi = (1-\epsilon)\pi_0 + \epsilon q; q \in \mathcal{Q}\}$, where $\pi_0$ is the original prior and $\mathcal{Q}$ represents any set of distributions to ensure that $\Gamma$ contains a plausible set of priors. Suppose the posterior distributions $\pi_0(\theta|x)$ and $q(\theta|x)$ exist (they do for a conjugate class $\mathcal{Q}$), and $m(x|\pi_0) = \int_\theta f(x|\theta)\pi_0(\theta)d\theta$ and $m(x|\pi)$ are the marginals obtained by integrating out $\theta$ with respect to the priors $\pi_0(\cdot)$ and $\pi(\cdot)$). It is straightforward to see that $\pi(\theta|x) = \lambda_{q,\epsilon}(x)\pi_0(\theta|x) + [1-\lambda_{q,\epsilon}(x)]q(\theta|x)$, where $\lambda_{q,\epsilon}(x) = \frac{(1-\epsilon)m(x|\pi_0)}{m(x|\pi)}$. It also follows that

$$(3.1) \qquad E_\pi(\theta|x) = \lambda_{q,\epsilon}(x)E_{\pi_0}(\theta|x) + [1-\lambda_{q,\epsilon}(x)]E_q(\theta|x).$$

However, posterior inference based on the MAP is affected, as the MAP involves the marginal likelihood instead of the conditional posterior. Let

$$MAP_{q,\epsilon}(A^*) = \frac{P_\pi(\boldsymbol{S} \mid A^*, \mathcal{M}_1)}{P_\pi(\boldsymbol{S} \mid \mathcal{M}_0)} = \frac{\int_{\mathcal{O}} P_\pi(\boldsymbol{S} \mid A^*, \mathcal{M}_1, \mathcal{O})\pi(\mathcal{O})d\mathcal{O}}{\int_{\mathcal{O}} P_\pi(\boldsymbol{S} \mid \mathcal{M}_0, \mathcal{O})\pi(\mathcal{O})d\mathcal{O}}.$$

If we are only interested in the effect of the prior $\pi(\Theta)$ for the motif column frequencies, the MAP for the contamination class $\Gamma$ reduces to

$$
\begin{aligned}
MAP_{q,\epsilon}(A^*) &= \frac{\int_{\mathcal{O}} P_\pi(\mathcal{S} \mid A^*, \mathcal{M}_1, \mathcal{O})\pi(\mathcal{O})d\mathcal{O}}{P(\mathcal{S} \mid \mathcal{M}_0)} \\
&= (1-\epsilon)MAP_{\pi_0}(A^*) + \epsilon MAP_q(A^*).
\end{aligned}
$$
(3.2)

### 3.1. Sensitivity analysis for Dirichlet prior for motif composition $\gamma$

Let us denote $\boldsymbol{\delta} = (\delta_A, \ldots, \delta_T)$ as the pseudocount vector for the Dirichlet distribution $q(\boldsymbol{\gamma})$. Also, let $\delta^* = \delta_A + \delta_T$. For simplicity, we assume that $\delta_A = \delta_T = \frac{\delta^*}{2}$, $\delta_C = \delta_G = 1 - \frac{\delta^*}{2}$. Then by varying $0 < \delta^* < 1$, we get a range of *contamination* distributions $q$. For $\epsilon \in (0,1)$, we evaluate the robustness of prior $\pi_0$, through the variability in magnitude of the posterior criteria in (3.1) and (3.2) over $\Delta = \{\boldsymbol{\delta}; 0 < \delta^* < 1\}$. Since we do not want the prior to dominate the data, we restrict the total pseudocount frequency $\sum_{j=1}^{4} \delta_j = 1$. We assume a single motif of width $w$ with observed frequency matrix $\mathcal{C}$. Let $\hat{\theta}_{ij} = E(\theta_{ij} \mid \mathcal{C})$. Then, (3.1) gives

$$
\begin{aligned}
\hat{\theta}_{ij} &= \lambda_q(x)\frac{c_{ij} + \gamma_{ij}}{\sum_{k=1}^{4}(c_{ik} + \gamma_{ik})} \\
&\quad + [1 - \lambda_q(x)]\frac{c_{ij} + \delta_{ij}}{\sum_{k=1}^{4}(c_{ik} + \delta_{ik})} \qquad 1 \leq i \leq w;\ 1 \leq j \leq 4. \\
&= \frac{1}{\sum_{j=1}^{4} c_{ij} + 1}\left[c_{ij} + \gamma_j + (\delta_j - \gamma_j)\left\{1 + \frac{1-\epsilon}{\epsilon}\prod_{i=1}^{w}\prod_{j=1}^{4}\left(\frac{\Gamma(c_{ij} + \gamma_j)}{\Gamma(c_{ij} + \delta_j)}\right)\right\}\right].
\end{aligned}
$$

Comparative studies were done separately to test the sensitivity of the (1) Dirichlet prior parameters $\boldsymbol{\beta}$ for word frequency and (2) product Dirichlet prior parameters $\boldsymbol{\gamma}$ for the column frequencies, as $\boldsymbol{\theta}$ and $\boldsymbol{\Theta}$ are posterior independent. Nucleotide composition in biological sequences are often asymmetric, e.g. a relatively higher frequency of C/G is observed for higher organisms. With known motif PWMs we study the effect on posterior inference when we use priors with pseudocount frequencies reflecting (1) uniform letter frequency (2) data composition-dependent letter frequencies or (3) a mixture with components having different letter frequencies.

Since we are mainly interested in whether the *consensus* motif sequence is affected by the choice of prior we study the variation of the highest frequency in each column $\hat{\theta}_i^* = \max_{1 \leq j \leq 4}\{\hat{\theta}_{ij}\}$. Corresponding to each choice of $\boldsymbol{\delta}$, we have a vector of *maximal* frequencies $\hat{\theta}^{*\delta} = (\hat{\theta}_1^{*\delta}, \ldots, \hat{\theta}_w^{*\delta})$. Let $\bar{\theta}_i^{*\delta} = \frac{1}{|\Delta|}\sum_\delta \hat{\theta}_i^{*\delta}$ denote the mean frequency for column $i$ ($i = 1, \ldots, w$). To summarize this information over all columns we calculate three distance-based measures, for each $\delta \in \Delta$ and $\epsilon \in (0,1)$–(i) Mean squared distance: $D_M^\delta = \frac{1}{w}\sum_{i=1}^{w}\left[\hat{\theta}_i^{*\delta} - \bar{\theta}_i^{*\delta}\right]^2$, (ii) Kullback-Leibler distance: $D_K^\delta = \sum_{i=1}^{w} \hat{\theta}_i^{*\delta} \log \frac{\hat{\theta}_i^{*\delta}}{\bar{\theta}_i^{*\delta}}$, and (iii) Entropy or Information content: $D_E^\delta = -\sum_{i=1}^{w} \hat{\theta}_i^{*\delta} \log \hat{\theta}_i^{*\delta}$.

All three measures $D_M, D_K$ and $D_E$ exhibit little variability, only after the third decimal place (plots not shown), slightly increasing with the contamination rate $\epsilon$. However, the MAP scores, calculated for the equal, data-dependent, 3-component and 9-component mixture prior indicate that using the mixture prior distribution makes the MAP most robust to the choice of prior (Figure 3). So if the results are to be evaluated using the MAP score at the mode, it seems to be safest to use a prior of Dirichlet mixtures instead of a single density. The most dramatic difference is observed in the GAL4 and MEF2 data on using a mixture Dirichlet prior– these are also the data sets having the lowest compositional similarity between the motif nucleotide frequency and background nucleotide frequency (Kullback-Leibler distance in the last column of Table 1). It is probable that the usual DA and the Gibbs sampler tend to pick up motifs that are closer in composition to the background for a similar reason. Using a mixture Dirichlet prior may be more sensitive towards detecting motifs that are vary compositionally from the background.

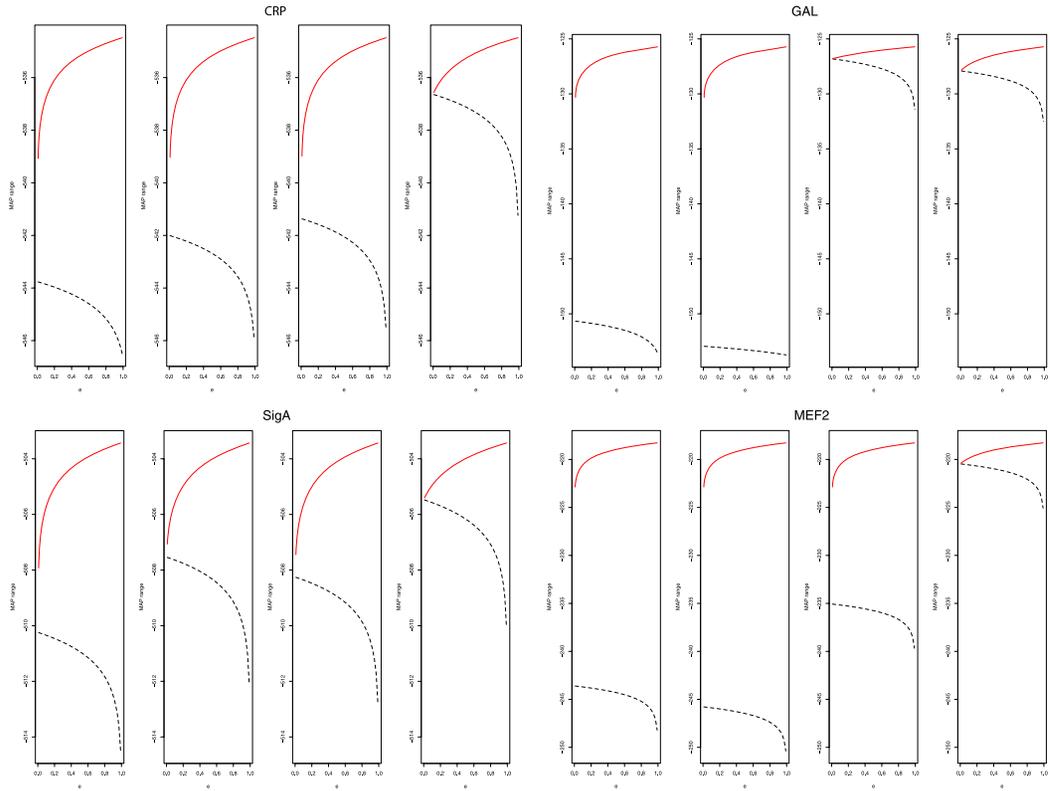

Fig 3. *Range in variability of MAP score (in log scale on vertical axis, may exclude an additive constant) at different levels of prior contamination, for data sets (1) CRP (2) GAL4 (3) $\sigma_A$ (TATAAT motif) (4) MEF2 for human skeletal muscle data set. The four panels for each data set correspond to taking the initial prior (i) Dirichlet with symmetric base pseudocounts (ii) Dirichlet with pseudocounts proportional to frequency in data (iii) 3-component mixture Dirichlet (iv) 9 component mixture Dirichlet. The horizontal axis on each plot represents the contamination proportion $\epsilon$.*

TABLE 1
*Motif and background letter frequencies for the 4 data sets*

| TF | Background | | | | Motif | | | | KL distance |
|---|---|---|---|---|---|---|---|---|---|
| | A | C | G | T | A | C | G | T | |
| CRP | 0.3021 | 0.1825 | 0.2090 | 0.3063 | 0.2841 | 0.1799 | 0.2140 | 0.3220 | 0.0011 |
| GAL4 | 0.3116 | 0.1914 | 0.1761 | 0.3208 | 0.1218 | 0.3908 | 0.2983 | 0.1891 | 0.2218 |
| $\sigma_A$ (TATAAT) | 0.3511 | 0.1465 | 0.1781 | 0.3243 | 0.4583 | 0.0699 | 0.0343 | 0.4375 | 0.1449 |
| MEF2 | 0.2205 | 0.2801 | 0.2715 | 0.2278 | 0.6047 | 0.0174 | 0.0262 | 0.3517 | 0.6531 |

## 3.2. Local sensitivity analysis

One problem in doing a careful robustness study in a high-dimensional problem is the limitation to mainly one type of high-dimensional prior distribution (the conjugate prior) for which analytical calculations are feasible. An alternative approach to judge sensitivity is to look for prior inputs that sharply change the posterior answers. Such *local* sensitivity analysis methods have been extensively developed and studied by Leamer [11] and Polasek [14]. In our framework, we observe that the posterior means are locally insensitive. Let $\mu_j = E(\rho_j \mid \mathcal{S}, \mathbf{A}) = \frac{N_j + \beta_j}{\sum_{k=1}^{D}(n_k + \beta_k)}$. Then, $\frac{\partial \mu_j}{\partial \beta_i} = \frac{\sum_{k \neq i}(n_k + \beta_k)}{\left[\sum_{k=1}^{D}(n_k + \beta_k)\right]^2} \leq 1$. A similar result holds for posterior means of $\boldsymbol{\Theta}$. However, parameter specifications for both $\boldsymbol{\gamma}$ and $\boldsymbol{\beta}$ may have a local impact on the MAP. Excluding constant terms, and using Stirling's approximation [1] to expand and simplify the $\Gamma$-functions, we have

$$\frac{\partial \log \text{MAP}}{\partial \gamma_j}$$
$$\approx \sum_i \left[ \log \frac{c_{ij} + \gamma_j}{\sum_k (c_{ik} + \gamma_k)} - \frac{1}{2} \left\{ \frac{1}{c_{ij} + \gamma_j} - \frac{1}{\sum_k (c_{ik} + \gamma_k)} \right\} \right]$$
$$+ w \left[ \log \frac{\sum_k \gamma_k}{\gamma_j} - \frac{1}{2} \left( \frac{1}{\sum_k \gamma_k} - \frac{1}{\gamma_j} \right) \right].$$

So the influence of $\gamma_j$ can be unbounded only for the second term– it increases with $w$ and as the ratio of $\sum_{k \neq j} \gamma_k$ to $\gamma_j$ increases. If $\sum_k \gamma_k = 1$, the second term is $w \left( \log \frac{1}{\gamma_j} + \frac{1}{2\gamma_j} - \frac{1}{2} \right)$, which can be made arbitrarily large if $\gamma_j$ is small. Under the restriction that $\sum_k \gamma_k = 1$, we will get most robust results if the components of $\boldsymbol{\gamma}$ are approximately equal.

We may also do a similar analysis for the local influence of the prior parameters for $\boldsymbol{\rho}$, i.e. for $\boldsymbol{\beta}$. Again, we assume there is only 1 motif type under the motif model

(index is $D$). If we differentiate logMAP with respect to $\beta_k$, $k < D$, we get

$$\log\left(\frac{N_{1k}}{N_{0k}}\right) + \frac{1}{2}\left[\frac{N_{1k} - N_{0k}}{(N_{1k} + \beta_k)(N_{0k} + \beta_k)}\right]$$

$$+ \log\left[\frac{\sum_{j=1}^{D-1}(N_{0j} + \beta_j)}{\sum_{j=1}^{D}(N_{1j} + \beta_j)}\right] + \log\left[\frac{\sum_{j=1}^{D}\beta_j}{\sum_{j=1}^{D-1}\beta_j}\right]$$

$$+ \frac{1}{2}\frac{\sum_{j=1}^{D-1} N_{0j} - \sum_{j=1}^{D} N_{1j} - \beta_D}{\left[\sum_{j=1}^{D}(N_{1j} + \beta_j)\right]\left[\sum_{j=1}^{D-1}(N_{0j} + \beta_j)\right]}$$

(3.3)
$$+ \frac{1}{2}\left(\frac{1}{\sum_{j=1}^{D-1}\beta_j} - \frac{1}{\sum_{j=1}^{D}\beta_j}\right),$$

where $\mathbf{N}_1$ and $\mathbf{N}_0$ are the frequencies under the motif ($\mathcal{M}_1$) and null ($\mathcal{M}_0$) models. Let $L$ denote the total length of the data. Then, $L = \sum_{j=1}^{D-1} N_{0j} = \sum_{j=1}^{D-1} N_{1j} + wN_D$, and $\sum_{j=1}^{D} N_{1j} = L - N_D(w-1)$. If $\beta_D < \sum_{j=1}^{D-1}\beta_j$, the influential part is:

$$\frac{1}{2}\frac{N_D(w-1) - \beta_D}{\left(L + \sum_{j=1}^{D-1}\beta_j - N_D(w-1) - \beta_D\right)\left(L + \sum_{j=1}^{D-1}\beta_j\right)}$$

$$- \log\left[1 - \frac{N_D(w-1) - \beta_D}{L + \sum_{j=1}^{D-1}\beta_j}\right].$$

Now, $N_D(w-1) < L$. So if $\beta_D < \sum_{j=1}^{D-1}\beta_j$, the influence of $\beta_k$ is seen to be negligible for $k < D$. On the other hand, if $\beta_D > \sum_{j=1}^{D-1}\beta_j$, then the influence of the third, fourth and fifth terms in (3.3) may be unbounded, but this is an unrealistic parametrization in practice, as the motif proportion is usually low compared to the size of the data. Finally, taking derivatives of logMAP with respect to $\beta_D$, we have,

$$\log\left[\frac{N_D + \beta_D}{\sum_{j=1}^{D}(N_{1j} + \beta_j)}\right] - \frac{\frac{1}{2}\sum_{j=1}^{D-1}(N_{1j} + \beta_j)}{(N_D + \beta_D)\sum_{j=1}^{D}(N_{1j} + \beta_j)}$$

$$+ \log\left(1 + \frac{\sum_{j=1}^{D-1}\beta_j}{\beta_D}\right) + \frac{1}{2}\left(\frac{1}{\beta_D} - \frac{1}{\sum_{j=1}^{D}\beta_j}\right).$$

The only influential term is the third, and the influence of $\beta_D$ can be unbounded if $\sum_{j=1}^{D-1}\beta_j \gg \beta_D$. So ideally, we want $\beta_D < \sum_{j=1}^{D-1}\beta_j$, but also that the difference is not extremely small. In other words, results are most robust to the choice of $\beta_D$ when motif proportions are comparatively large. Having prior knowledge of the true motif frequency would make it easier to get more accurate results– by choosing $\beta_D$ to give the expected motif frequency close to the prior knowledge. In practice, the choice of $\beta_D$ is seen to have a great impact on motif detection and evaluation, and this remains one of the greatest stumbling blocks of this probabilistic framework.

## 4. Discussion and conclusions

A MAP criterion is proposed for model selection in the motif discovery problem. Analytical and computational investigations establish conditions for the the MAP criterion to asymptotically predict the correct number of motifs and attaining these conditions is seen to be feasible in a real scenario. We investigate the sensitivity of the MAP criterion to prior specification and provide guidelines for choosing these hyper-parameters to maximize robustness. The MAP is seen to perform well as a model selection criterion in preliminary studies, paving the way for further analysis of its performance when adapted to more complex realistic models.

## Appendix A: Outline of proof for Theorem 2.1

From expression (2.3), we have

$$
\begin{aligned}
MAP(A^0) &\approx \frac{\prod_{i=1}^{d}(N_{1i}+\beta_i-1)^{N_{1i}+\beta_i-\frac{1}{2}}(N_{1,d+1}+\beta_{d+1}-1)^{N_{1,d+1}+\beta_{d+1}-\frac{1}{2}}}{\prod_{i=1}^{d}(N_{0i}+\alpha_i-1)^{N_{0i}+\alpha_i-\frac{1}{2}}} \\
&\quad \times \frac{\left(\sum_{i=1}^{d}[N_{0i}+\alpha_i]-1\right)^{\sum_{i=1}^{d}[N_{0i}+\alpha_i]-\frac{1}{2}}}{\left(\sum_{i=1}^{d+1}[N_{1i}+\beta_i]-1\right)^{\sum_{i=1}^{d+1}[N_{1i}+\beta_i]-\frac{1}{2}}} \\
&\quad \times \prod_{i=1}^{w}\frac{\prod_{j=1}^{d}\Gamma(c_{ij}+\gamma_j)}{\Gamma\left(\sum_{j=1}^{d}[c_{ij}+\gamma_j]\right)} \times k(\alpha,\beta,\gamma),
\end{aligned}
$$

where

$$
k(\alpha,\beta,\gamma) = \frac{\left(\sum_{i=1}^{d+1}\beta_i-1\right)^{\sum_{i=1}^{d+1}\beta_i-\frac{1}{2}}}{\left(\sum_{i=1}^{d}\alpha_i-1\right)^{\sum_{i=1}^{d}\alpha_i-\frac{1}{2}}} \times \frac{\prod_{i=1}^{d}(\alpha_i-1)^{\alpha_i-\frac{1}{2}}}{\prod_{i=1}^{d+1}(\beta_i-1)^{\beta_i-\frac{1}{2}}} \times \frac{\Gamma\left(\sum_{j=1}^{d}\gamma_j\right)}{\prod_{i=1}^{d}\Gamma(\gamma_j)}.
$$

Denote $\sum_{j=1}^{d}c_{ij} = N_{1,d+1} = m$, $\sum_{i=1}^{d}N_{0i} = \sum_{i=1}^{d}N_{1i} + mw = N \Rightarrow \sum_{i=1}^{d+1}N_{1i} = N - (w-1)m$. Then,

$$
\begin{aligned}
MAP(A^0) &\approx k(\alpha,\beta,\gamma) \times \prod_{i=1}^{d}\frac{\left(1+\frac{\beta_i-1}{N_{1i}}\right)^{N_{1i}+\beta_i-\frac{1}{2}}}{\left(1+\frac{\alpha_i-1}{N_{0i}}\right)^{N_{0i}+\alpha_i-\frac{1}{2}}} \\
&\quad \times (m+\beta_{d+1}-1)^{m+\beta_{d+1}-\frac{1}{2}}\prod_{i=1}^{d}\frac{N_{1i}^{N_{1i}+\beta_i-\frac{1}{2}}}{N_{0i}^{N_{0i}+\alpha_i-\frac{1}{2}}} \\
&\quad \times \left\{\left(1+\frac{\sum_{i=1}^{d}\alpha_i-1}{N}\right)^{N+\sum_{i=1}^{d}\alpha_i-\frac{1}{2}} N^{N+\sum_{i=1}^{d}\alpha_i-\frac{1}{2}} \right. \\
&\quad\quad \times \left(1+\frac{\sum_{i=1}^{d+1}\beta_i-1-(w-1)m}{N}\right)^{N-(w-1)m+\sum_{i=1}^{d+1}\beta_i-\frac{1}{2}} \\
&\quad\quad \left. \times N^{N-(w-1)m+\sum_{i=1}^{d+1}\beta_i-\frac{1}{2}}\right\} \\
\text{(A.1)} &\quad \times \prod_{i=1}^{w}\frac{\prod_{j=1}^{d}(c_{ij}+\gamma_j-1)^{c_{ij}+\gamma_j-\frac{1}{2}}}{\left(m+\sum_{j=1}^{d}\gamma_j-1\right)^{m+\sum_{j=1}^{d}\gamma_j-\frac{1}{2}}}.
\end{aligned}
$$

Let the symbol $\cong$ denote "equal in order of convergence", i.e. ignoring terms that converge to a constant in the limit, as $N \to \infty, N_{0i} \to \infty, N_{1i} \to \infty$. Also, as $m \to \infty$, $\frac{m}{N} \xrightarrow{P} c$. Then, taking logarithms, expression (A.1)

$$\cong \sum_{i=1}^{d}(\beta_i - \alpha_i) + \sum_{i=1}^{d} \alpha_i - 1$$

$$+ \left(\sum_{i=1}^{d+1} \beta_i - \frac{1}{2}\right) \log[1 - c(w-1)] + \left[(w-1)m + \sum_{i=1}^{d} \alpha_i - \sum_{i=1}^{d+1} \beta_i\right] \log N$$

$$- \frac{w(d-1)}{2} \log\left(1 + \frac{\sum_{j=1}^{d} \gamma_j}{m}\right) + \left(m + \beta_{d+1} - \frac{w(d-1)+1}{2}\right) \log m$$

$$+ \sum_{i=1}^{d} \left[\log N_{1i}^{N_{1i}+\beta_i - \frac{1}{2}} - \log N_{0i}^{N_{0i}+\alpha_i - \frac{1}{2}}\right]$$

$$+ e^{\beta_{d+1}} - 1 - [1 - c(w-1)]N \log[1 - c(w-1)]$$

(A.2) $$+ \sum_{i=1}^{w} \sum_{j=1}^{d} \left(c_{ij} + \gamma_j - \frac{1}{2}\right) \log\left(\frac{c_{ij} + \gamma_j - 1}{m + \sum_{j=1}^{d} \gamma_j - 1}\right),$$

since

$$\left(1 - \frac{(w-1)m + \sum_{i=1}^{d+1} \beta_i - \frac{1}{2}}{N}\right)^N \xrightarrow{P} [1 - c(w-1)]^N \to 0$$

(same for limit in $m$). Under the conditions of the statement, $N_{0i} \to \infty$, $N_{1i} \to \infty$, $\frac{N_{0i} - N_{1i}}{m} \xrightarrow{P} k_i w$, $\frac{N_{0i}}{N} \xrightarrow{P} \theta_{0i}$. Now, assume that the motif is "exact", i.e. $c_{ij} = m$ for some $j$, every $i$. Now let us denote $\gamma^+ = \max_{1 \leq j \leq d} \gamma_j$, $\gamma^- = \min_{1 \leq j \leq d} \gamma_j$. Then, the last term of (A.2) is greater than

(A.3) $$w\left(m + \gamma^- - \frac{1}{2}\right) \log(m + \gamma^- - 1) \cong Nwc \log c + w\left(cN + \gamma^- - \frac{1}{2}\right) \log N.$$

Taking the limits of all non-constant terms in expression (A.2), and simplifying, we then have

$$\log MAP(A^0)$$

$$\cong \left[\sum_{i=1}^{d}(\theta_{0i} - k_i cw) \log(\theta_{0i} - k_i cw)\right.$$

$$\left. - \sum_{i=1}^{d} \theta_{0i} \log \theta_{0i} + c \log c - [1 - c(w-1)] \log[1 - c(w-1)]\right] N$$

(A.4) $$+ \left[w\gamma^- - w\sum_{j=1}^{d} \gamma_j - \frac{1}{2}\right] \log N.$$

So a sufficient condition for $MAP(A^0) \xrightarrow{P} \infty$ at an exponential rate is

$$c \log c + \sum_{i=1}^{d}(\theta_{0i} - k_i cw) \log(\theta_{0i} - k_i cw)$$

(A.5) $$- \sum_{i=1}^{d} \theta_{0i} \log \theta_{0i} - [1 - c(w-1)] \log[1 - c(w-1)] > 0.$$

Now we will show that (A.5) actually is a necessary condition. Assume that condition (A.5) is not satisfied. If the first term in (A.4) is $< 0$, then $\log MAP(A^0) \to -\infty$. If the first term is zero, using (A.3), we see that,

$$
\begin{aligned}
\log MAP(A^0) &\geq \left[ w\gamma^- - w \sum_{j=1}^{d} \gamma_j - \frac{1}{2} \right] \log N, \text{ and} \\
\log MAP(A^0) &\leq \left[ w\gamma^+ - w \sum_{j=1}^{d} \gamma_j - \frac{1}{2} \right] \log N.
\end{aligned}
\tag{A.6}
$$

In (A.6), the coefficient of $\log N$ is negative. Thus, as $N \to \infty$, $MAP(A^0)$ is bounded by 0. So condition (A.5) is also a *necessary* condition for $MAP(A^0) \to \infty$.

## Appendix B: Outline of proof for Theorem 2.2

Using the Lagrangian method with the restrictions ($\sum_{i=1}^{d} \theta_{0i} = 1$) and ($\sum_{i=1}^{d} k_i = 1$), we get $\theta_{0i} = k_i = \frac{1}{d}$ as the optimal solution. To verify that this is the *maximum* value, let

$$F(\boldsymbol{\theta}_0, \boldsymbol{k}) = \sum_{i=1}^{d} (\theta_{0i} - k_i cw) \log(\theta_{0i} - k_i cw) - \sum_{i=1}^{d} \theta_{0i} \log \theta_{0i}.$$

Then, $\frac{\partial^2 F}{\partial \theta_{0i}^2} = \frac{1}{\theta_{0i} - k_i cw} - \frac{1}{\theta_{0i}} > 0$, $\frac{\partial^2 F}{\partial k_i^2} = \frac{c^2 w^2}{\theta_{0i} - k_i cw} > 0$, and $\frac{\partial^2 F}{\partial \theta_{0i} \partial k_i} = \frac{-cw}{\theta_{0i} - k_i cw} < 0$. The determinant of the Hessian matrix is $\|H\| = \left\| \frac{\partial^2 F}{\partial \boldsymbol{a} \partial \boldsymbol{k}} \right\| = \left\| \begin{array}{cc} A & B \\ B & K \end{array} \right\|$, where $A = diag(A_i)$; $(A_i = \frac{\partial^2 F}{\partial \theta_{0i}^2})$, $B = diag(B_i)$; $(B_i = \frac{\partial^2 F}{\partial \theta_{0i} \partial k_i})$ and $K = diag(K_i)$; $(K_i = \frac{\partial^2 F}{\partial k_i^2})$. Hence $\|H\| = |A||K - BA^{-1}B| = |A||diag(M_i)|$, where

$$
\begin{aligned}
M_i &= \frac{c^2 w^2}{\theta_{0i} - k_i cw} - \frac{c^2 w^2}{(\theta_{0i} - k_i cw)^2} \times \frac{1}{\left( \frac{1}{\theta_{0i} - k_i cw} - \frac{1}{\theta_{0i}} \right)} \\
&= \frac{cw}{\theta_{0i} - k_i cw} \left( cw - \frac{\theta_{0i}}{k_i} \right) < 0.
\end{aligned}
$$

Hence, $\|H\|$ is negative definite, and the expression in (2.5) corresponds to the maximum attained by $F$.

## Appendix C: Outline of proof for Theorem 2.3

Let $n$ denote the sum of word frequencies under model $\mathcal{M}_q$, i.e. $n = \sum_{i=1}^{d} N_{1i}$. Let $m$ be the number of word instances of type $q+1$ and $l$ be the corresponding number of type $q$, $q \geq d$ where the alphabet size is $d$ as previously. Ideally, if model $\mathcal{M}_{q+1}$ is true, when $n \to \infty$, with $\frac{N_{q+1}}{n} \to \rho_{q+1}$, $\log MAP_{q+1} - \log MAP_q \to \infty$; while

$\log MAP_{q+1} - \log MAP_q \to 0$ (or is bounded above) when $\frac{N_{q+1}}{n} = 0$. By definition,

$$\log MAP^m_{q+1} - \log MAP^l_q$$

$$= \log \left[ \frac{\prod_{i=1}^{q+1} \Gamma(N_{2i} + \beta_i) \Gamma\left(\sum_{i=1}^{q+1} \beta_i\right)}{\Gamma\left(\sum_{i=1}^{q+1}[N_{2i} + \beta_i]\right) \prod_{i=1}^{q+1} \Gamma(\beta_i)} \right]$$

$$- \log \left[ \frac{\prod_{i=1}^{q} \Gamma(N_{1i} + \beta_i) \Gamma\left(\sum_{i=1}^{q} \beta_i\right)}{\Gamma\left(\sum_{i=1}^{q}[N_{1i} + \beta_i]\right) \prod_{i=1}^{q} \Gamma(\beta_i)} \right]$$

$$+ \log \left[ \prod_{i=1}^{w} \frac{\prod_{j=1}^{d} \Gamma(c_{ijk} + \gamma_j)}{\Gamma\left(\sum_{j=1}^{d} \Gamma[c_{ijk} + \gamma_j]\right)} \right]$$

(C.1)
$$+ w \log \left[ \frac{\Gamma\left(\sum_{j=1}^{d} \gamma_j\right)}{\prod_{j=1}^{d} \Gamma(\gamma_j)} \right], \quad \text{(where } k = q + 1 - d\text{)}.$$

If $\frac{N_{2,q+1}}{N} = 0$, then in (C.1), $\frac{N_{1i}}{N} - \frac{N_{2i}}{N} \to 0$ for $i \leq q$. Then (C.1) becomes $\log[\Gamma(\sum_{i=1}^{q+1} \beta_i)] - \log[\Gamma(\sum_{i=1}^{q} \beta_i)] - \log \Gamma(\beta_{q+1}) = \log \text{Beta}(\sum_{i=1}^{q} \beta_i, \beta_{q+1})$. A sufficient condition for $\log MAP_{q+1} - \log MAP_q < 0$ is $\beta_{q+1} > 1$, $\sum_{i=1}^{q} \beta_i > 1$. Now we derive the monotonicity property of the MAP score, i.e. when $\frac{N_{q+1}}{N} \to \rho_{q+1}$, under model $\mathcal{M}_{q+1}$, $\log MAP_{q+1} - \log MAP_q \to \infty$. Leaving out the constant terms, (C.1) reduces to

$$\log \left[ \frac{\prod_{i=1}^{d} \Gamma(N_{2i} + \beta_i) \Gamma(N_{2,q+1} + \beta_{q+1}) \Gamma\left(\sum_{i=1}^{q}[N_{1i} + \beta_i]\right)}{\prod_{i=1}^{d} \Gamma(N_{1i} + \beta_i) \Gamma\left(\sum_{i=1}^{q+1}[N_{2i} + \beta_i]\right)} \right]$$

$$+ \sum_{i=1}^{w} \log \left[ \frac{\prod_{j=1}^{d} \Gamma(c_{ijk} + \gamma_j)}{\Gamma\left(\sum_{j=1}^{d} c_{ijk} + \sum_{j=1}^{d} \gamma_j\right)} \right].$$

Let us denote $m = N_{2,q+1} = \sum_{j=1}^{d} c_{ijk} \forall i = 1, \ldots, w$. Then $\sum_{i=1}^{q} N_{2i} + mw = \sum_{i=1}^{q} N_{1i} = n$. Using Stirling's approximation, (C.1) reduces to

$$\sum_{i=1}^{d} \left[ \left(N_{2i} + \beta_i - \frac{1}{2}\right) \log(N_{2i} + \beta_i - 1) - \left(N_{1i} + \beta_i - \frac{1}{2}\right) \log(N_{1i} + \beta_i - 1) \right]$$

$$- \left(m + \beta_{q+1} - \frac{1}{2}\right) \log(m + \beta_{q+1} - 1)$$

$$+ \left(\sum_{i=1}^{q}[N_{1i} + \beta_i] - \frac{1}{2}\right) \log \left(\sum_{i=1}^{q}[N_{1i} + \beta_i] - 1\right)$$

$$- \left(\sum_{i=1}^{q+1}[N_{2i} + \beta_i] - \frac{1}{2}\right) \log \left(\sum_{i=1}^{q+1}[N_{2i} + \beta_i] - 1\right)$$

$$+ \sum_{i=1}^{w} \left[ \sum_{j=1}^{d} \left(c_{ijk} + \gamma_j - \frac{1}{2}\right) \log(c_{ijk} + \gamma_j - 1) \right.$$

(C.2)
$$\left. - \left(m + \sum_{j=1}^{d} \gamma_j - \frac{1}{2}\right) \log \left(m + \sum_{j=1}^{d} \gamma_j - 1\right) \right],$$

since, for $i = d+1, \ldots, q$, we have $N_{1i} = N_{2i}$. Leaving out constants and terms tending towards 0 as $N_{1i} \to \infty$, $N_{2i} \to \infty$, we can write (C.2) as

$$\sum_{i=1}^{d} \left[\left(N_{2i} + \beta_i - \frac{1}{2}\right) \log(N_{2i}) - \left(N_{1i} + \beta_i - \frac{1}{2}\right) \log(N_{1i})\right]$$

$$+ \left(m + \beta_{q+1} - \frac{1}{2}\right) \log(m + \beta_{q+1} - 1)$$

$$+ \left(n + \sum_{i=1}^{q} \beta_i - \frac{1}{2}\right) \log\left(1 + \frac{\sum_{i=1}^{q} \beta_i - 1}{n}\right)$$

$$- \left[n - (w-1)m + \sum_{i=1}^{q+1} \beta_i - \frac{1}{2}\right] \log\left[1 - \frac{(w-1)m + \sum_{i=1}^{q+1} \beta_i + 1}{n}\right]$$

$$+ [m(w-1) - \beta_{q+1}] \log(n)$$

$$+ \sum_{i=1}^{w} \left[\sum_{j=1}^{d} \left(c_{ijk} + \gamma_j - \frac{1}{2}\right) \log(c_{ijk} + \gamma_j - 1)\right.$$

(C.3) $$\left. - \left(m + \sum_{j=1}^{d} \gamma_j - \frac{1}{2}\right) \log\left(m + \sum_{j=1}^{d} \gamma_j - 1\right)\right].$$

Now we evaluate each of the terms as $n \to \infty, m \to \infty$, with $\frac{m}{n} \xrightarrow{P} \rho_{q+1}$. Assume that the proportion contributed to motif $(q+1)$ by each base $i$, $i = 1, \ldots, d$ is $\kappa_i$, such that $\frac{N_{1i} - N_{2i}}{m} \xrightarrow{P} w\kappa_i$, $i = 1, \ldots, d$, $\sum_{i=1}^{d} \kappa_i = 1$ and, as before, $n - mw = \sum_{i=1}^{q} N_{2i}$. Also, assume that all instances of the motif are exact, i.e. $c_{ijk} = m$ for some $j \in 1, \ldots, d$ $\forall i = 1, \ldots, w$. Also, for all $i, j$, let $\gamma(i) = \gamma_j$ if $c_{ijk} = m$. So, finally, simplifying terms in (C.3), we get,

$$n\left[\sum_{i=1}^{d} \{(\rho_i - \kappa_i w \rho_{q+1}) \log(\rho_i - \kappa_i w \rho_{q+1}) - \rho_i \log \rho_i\}\right.$$

$$\left. + [1 - (w-1)\rho_{q+1}] \log[1 - (w-1)\rho_{q+1}] + \rho_{q+1} \log \rho_{q+1}\right]$$

(C.4) $$+ n \log n[-w + 1 + w - 1]\rho_{q+1} + \log n \left[\sum_{i=1}^{w} \gamma(i) - w \sum_{j=1}^{d} \gamma_j\right].$$

Let $r$ denote

$$\sum_{i=1}^{d} \{(\rho_i - \kappa_i w \rho_{q+1}) \log(\rho_i - \kappa_i w \rho_{q+1}) - \rho_i \log \rho_i\}$$

$$+ \rho_{q+1} \log \rho_{q+1} + [1 - (w-1)\rho_{q+1}] \log[1 - (w-1)\rho_{q+1}].$$

Then (C.4) $\xrightarrow{P} \infty$ at a rate exponential in $rn$ iff $r > 0$. If $r = 0$, the coefficient of $\log n$, $\sum_{i=1}^{w} \gamma(i) - w \sum_{j=1}^{d} \gamma_j \leq 0$ so $\log MAP_{q+1} - \log MAP_q \xrightarrow{P} 0$. Also, by a similar argument as in Theorem 2.2, the maximum value is attained when $\kappa_i = \rho_i = \frac{1}{d}$, for $i = 1, \ldots, d$.

**Acknowledgments.** The author is grateful to Jun S. Liu for helpful comments.